\documentclass{article}
\usepackage[latin1]{inputenc}
\usepackage{ngerman}
\usepackage{url}
\usepackage{amssymb}
\usepackage[longtable]{ednotes}
\setpagewiselinenumbers
\modulolinenumbers[5]
\firstlinenumber{5}
\usepackage{fancyhdr}
\def\seitennummer{\thepage}

\usepackage{scrtime}
\makeatletter
\long\def\wennleer#1#2#3{\def\@tempa{#1}\ifx\@tempa\empty#2\else#3\fi}
\makeatother
\usepackage{natbib}
 \let\cite\citep
 
  \bibpunct{(}{)}{;}{}{}{,~}
\let\Bibitem\bibitem
\renewcommand*{\bibitem}[2][]{\def\diesesigle{#2}\Bibitem[#1]{#2}}
\def\letzterautor{}
\def\fontsaufbibfonts{\def\bibnamefont##1{\textsc{##1}}\def\bibfnamefont##1{\textsc{##1}}}
\def\autor#1{\def\letzterautor{#1}\IfFirstOnPageThenElse{\letzterautor}{{\fontsaufbibfonts #1}}{---} (\citeyear{\diesesigle})}
\def\drs{\IfFirstOnPageThenElse{\letzterautor}{{\fontsaufbibfonts\letzterautor}}{---}~(\citeyear{\diesesigle})}
\def\drsit{\IfFirstOnPageThenElse{\letzterautor}{{\fontsaufbibfonts\letzterautor}}{---}\ }
\def\autorit#1{\def\letzterautor{#1}\IfFirstOnPageThenElse{\letzterautor}{{\name{#1}}}{---}}
\makeatletter 
\newcount\c@GP@something
\def\GPbadtrue {\global\let\ifGPbad\iftrue}  
\def\GPbadfalse{\global\let\ifGPbad\iffalse}
\def\@GP@occname{@GP@no\number\c@GP@something} 
\def\@GP@mark#1{%
  \global\advance\c@GP@something\@ne 
  \global \expandafter \mathchardef 
    \csname \@GP@occname \endcsname #1\relax
}
\newcommand\GuessPageNo{
  \global\c@GP@something\z@ 
  \gdef\GuessPageNo{%
    \global\advance\c@GP@something\@ne 
    \protected@write\@auxout{%
      \let\number\relax}{%
        \string\@GP@mark{\number\c@page}}%
    \@ifundefined\@GP@occname{%
      \GPbadtrue \G@refundefinedtrue}{%
      \GPbadfalse 
      \global \expandafter \let \expandafter 
        \c@GuessedPageNo \csname\@GP@occname\endcsname 
    }%
  }%
  \GuessPageNo 
}
\AtEndDocument{%
  \global\c@GP@something\z@ 
  \def\@GP@mark#1{%
    \global\advance\c@GP@something\@ne 
    \expandafter \ifx \csname \@GP@occname \endcsname \relax 
    \else 
      \expandafter \ifnum 
        \csname \@GP@occname \endcsname=#1\relax 
      \else 
        \@tempswatrue 
      \fi 
    \fi 
  }%
}
\def\@GP@catname#1{n@GP@cat@#1}
\newcommand\IfFirstOnPageThenElse[1]{%
  \GuessPageNo 
  \edef\@tempa{n@GP@cat#1}%
  \@ifundefined\@tempa{%
    \GPbadtrue}{%
    \ifnum \csname\@tempa\endcsname < \c@GuessedPageNo 
      \GPbadtrue
    \fi 
  }%
  \ifGPbad
    \global\expandafter\mathchardef 
      \csname \@tempa \endcsname \c@GuessedPageNo 
    \expandafter\@firstoftwo 
  \else 
    \expandafter\@secondoftwo 
  \fi 
} 
\makeatother 
\def\phi{\varphi}
\usepackage[]{paralist}

\long\def\sina#1#2{#1\footnote{\ensuremath{\mathfrak{S}} writes \anf{#2}.}}
\long\def\anf#1{\glqq #1\grqq}
\def\seiteo#1{$\mid$\marginpar[]{\small #1}}
\begin{document}
\thispagestyle{plain}
\noindent\texttt{MODIFIED DRAFT VERSION, 30.\,Nov.\,2016\\
Comments welcome!}\par\vspace{2em}
\begin{center}
\sc
Richard Dedekind: Brief an Keferstein / Letter to Keferstein \par
\it
edited by Christian Tapp
\end{center}\par\vspace{1em}
This is a critical edition of the famous letter, the mathematician 
Richard Dedekind (1831--1916) has sent to one of his critiques, namely H.~Keferstein on 27
February 1890.\par
The original German version was first published by M.-A.~Sinaceur in \cite{Sinaceur1974}.
The edition presented here relies on Sinaceur's one but corrects some 
inaccuracies.\par
For more biographical information on Dedekind see \cite{OConnorRobertson1998}. For an account of his foundational research 
see \cite{SiegSchlimm2005}.\par
The \emph{edition is designed} as follows. Abbreviations are expanded by square brackets
(Abb[reviation]). Textcritical remarks (i.e., relevant differences to the edition 
\cite{Sinaceur1974} [=\ensuremath{\mathfrak{S}}]) are given in footnotes. 
Dedekind's old style German orthography has been reproduced in the original
way in order to provide the reader with the \anf{original feeling} of 
Dedekind's writing. For example, he always writes \anf{ß} instead of \anf{ss}, 
uses \anf{c} instead of \anf{k} in words with Latin origin, and he
adheres to writing \anf{th} instead of \anf{t} like in \anf{Thal} (old) instead of 
\anf{Tal} (modern).\par\vspace{1em}
\begin{center} \rule{10em}{0.1pt} \end{center}\par\vspace{1em}
1890.\,2.\,27. Herrn Oberlehrer Dr.\,H.\,Keferstein. Hamburg\par\vspace{0.5ex}\noindent
\mbox{}\rule{\linewidth}{0.1pt}\par\vspace{1ex}
Hochgeehrter Herr Doctor!\par\noindent
Für Ihren \sina{freundlichen}{freudlichen} Brief vom 14.\,d[ieses] M[onats] und die darin ausgesprochene Bereitwilligkeit, meiner Entgegnung Gehör zu verschaffen, sage ich Ihnen meinen besten Dank. Doch möchte ich Sie bitten, in dieser Angelegenheit Nichts zu übereilen und erst dann einen Entschluß zu faßen, nachdem Sie, wenn Sie \sina{Zeit}{die Zeit} dazu haben, die wichtigsten Erklärungen und Beweise in meiner Zahlenschrift noch einmal genau gelesen und durchdacht haben. Ich halte es nämlich für höchst wahrscheinlich, daß Sie sich dann in allen Puncten zu meiner \sina{Auffaßung}{Anffassung} und Behandlung des Gegenstandes bekehren werden, und gerade hierauf würde ich bei weitem den größten Werth legen, weil ich überzeugt bin, daß Sie wirklich ein tiefes Intereße für die Sache hegen.\par
Um diese Annäherung wo möglich zu befördern, möchte ich Sie bitten, dem folgenden Gedankengange, der die Genesis meiner Schrift darstellt, Ihre Aufmerksamkeit zu schenken. Wie ist meine Schrift entstanden? Gewiß nicht in einem \sina{Tage}{Zuge}, sondern sie ist eine nach langer Arbeit aufgebaute Synthese, die sich auf eine \seiteo{2} vorausgehende Analyse der Reihe der natürlichen Zahlen stützt, so wie diese sich, gewißermaßen erfahrungsmäßig, unserer Betrachtung darbietet. Welches sind die von einander unabhängigen Grundeigenschaften dieser Reihe $\mathcal{N}$, d.\,h. diejenigen Eigenschaften, welche sich nicht \sina{aus einander}{auseinander} ableiten laßen, aus denen aber alle anderen \sina{folgen?}{folgen.} Und wie muß man diese Eigenschaften \sina{ihres spezifisch arithmetischen}{Ihres spezifischarithmetischen} Charakters entkleiden, der Art, daß sie sich \sina{allgemeineren}{allgemeinen} Begriffen und solchen Tätigkeiten des Verstandes unterordnen, \emph{ohne} welche überhaupt kein Denken möglich ist, \emph{mit} welchen aber auch die Grundlage gegeben ist für die Sicherheit und Vollständigkiet der Beweise, wie für die Bildung widerspruchsfreier Begriffserklärungen?\par
Stellt man die Frage in dieser Weise, so wird man, wie ich glaube, mit Gewalt auf folgende Thatsachen gedrängt:\par
\begin{asparaenum}[1)]
\item Die Zahlenreihe $\mathcal{N}$ ist ein \emph{System} von Individuen oder Elementen, die 
\sina{Zahlen}{\anf{Zahlen}} heißen. Dies führt zur allge\seiteo{3}meinen Betrachtung von Systemen überhaupt (§. 1
meiner Schrift).
\item Die Elemente des Systems $\mathcal{N}$ stehen in gewißer Beziehung \sina{zu einander}{zueinander}, es herrscht eine gewiße Ordnung, die zunächst darin besteht, daß zu jeder bestimmten Zahl $n$ eine bestimmte Zahl $n'$, die folgende oder nächst größere Zahl \emph{gehört}. \sina{Dies}{Dieses} führt zur Betrachtung des allgemeinen Begriffs einer \emph{Abbildung} $\phi$ eines Systems (§. 2), und da das Bild $\phi(n)$ einer jeden Zahl $n$ wieder eine \sina{\emph{Zahl}}{Zahl} $n'$, also $\phi(\mathcal{N})$ Theil von $\mathcal{N}$ ist, so handelt es sich hier um die Abbildung $\phi$ eines Systems $\mathcal{N}$ \emph{in sich selbst}, welche also allgemein zu untersuchen ist (§. 4).
\item Auf verschiedene Zahlen $a$, $b$ folgen auch verschiedene Zahlen $a'$, $b'$; die Abbildung $\phi$ \sina{hat}{\emph{hat}} also den Charakter der Deutlichkeit oder \emph{Ähnlichkeit} (§. 3).
\item Nicht jede Zahl ist eine folgende Zahl $n'$, d.\,h. $\phi(\mathcal{N})$ ist echter Theil von $\mathcal{N}$, und hierin besteht (in Verbindung mit dem Vorhergehenden) die \emph{Unendlichkeit} der Zahlenreihe $\mathcal{N}$ (§. 5).
\item Und zwar ist die Zahl $1$ die \emph{einzige} Zahl, welche sich nicht in $\phi(\mathcal{N})$ findet. Hiermit sind \seiteo{4} diejenigen Thatsachen aufgezählt, in welchen Sie (S.\,124,\,Z.\,8--14) den vollständigen Charakter eines geordneten einfach unendlichen Systems $\mathcal{N}$ erblicken.
\item Aber ich habe in meiner Entgegnung (III) gezeigt, daß diese Thatsachen noch lange nicht ausreichen, um das Wesen der Zahlenreihe $\mathcal{N}$ vollständig zu \sina{erfaßen. Alle}{erfassen; Alle} diese Thatsachen würden auch noch für jedes System $\mathcal{S}$ gelten, welches außer der Zahlenreihe $\mathcal{N}$ noch ein System $\mathcal{T}$ von beliebigen anderen Elementen $t$ enthält, auf \sina{welches}{welche} die Abbildung $\phi$ sich stets so ausdehnen läßt, daß sie den Charakter der Ähnlichkeit behält, und daß $\phi(\mathcal{T}) = \mathcal{T}$ wird. Aber ein solches System $\mathcal{S}$ ist offenbar etwas ganz Anderes, als unsere Zahlenreihe $\mathcal{N}$, und ich könnte es so wählen, daß in ihm kaum ein einziger der arithmetischen Sätze bestehen bliebe. Was muß also zu den bisherigen Thatsachen noch hinzu kommen um unser System $\mathcal{S}$ von solchen fremden, alle Ordnung störenden Eindringlingen $t$ wieder zu reinigen und auf $\mathcal{N}$ zu \sina{beschränken?}{beschränken.} \sina{Dies}{Dieses} war einer der schwierigsten Puncte meiner Analyse, und seine Überwindung hat ein langes Nachdenken erfordert. Setzt man die Kenntniß der natürlichen Zahlenreihe $\mathcal{N}$ schon voraus und \seiteo{5} erlaubt sich demgemäß eine arithmetische Ausdrucksweise, so hat man ja leichtes Spiel; man braucht nur zu sagen: ein Element $n$ gehört dann und nur dann der Reihe $\mathcal{N}$ an, wenn ich, von dem Element 1 ausgehend, durch beständig wiederholtes Weiterzählen, d.\,h. durch eine \sina{endliche}{unendliche} Anzahl von Wiederholungen der Abbildung $\phi$ (vergl. den Schluß von 131 meiner Schrift) wirklich einmal zu dem Element $n$ gelange, während ich auf diese Weise niemals zu einem der Reihe $\mathcal{N}$ fremden Element $t$ gelange. Aber dieser Weg, den Unterschied zwischen den aus $\mathcal{S}$ wieder auszutreibenden Elementen $t$ und den allein bei zu behaltenden Elementen $n$ zu charakterisiren, ist doch für unseren Zweck gänzlich unbrauchbar, er enthielte ja einen \textsl{circulus} \sina{\textsl{vitiosus}}{viciosus} der schlimmsten und auffälligsten Art. Die bloßen Worte \anf{endlich einmal hinkommen} tun es auch natürlich nicht, sie würden nicht mehr helfen als etwa die Worte \textsl{\anf{Karam sipo tatura}}, die ich augenblicklich erfinde, ohne ihnen einen deutlich erklärten Sinn zu geben. Also: wie kann ich, ohne irgend welche arithmetische \sina{Kenntniß}{Erkenntnis} vorauszusetzen, den Unterschied zwischen \seiteo{6} den Elementen $n$ und $t$ unfehlbar begrifflich bestimmen? Ganz allein durch die Betrachtung der \emph{Ketten} (37, 44 meiner Schrift), durch diese aber auch vollständig! Will ich meinen Kunstausdruck \anf{Kette} vermeiden, so werde ich sagen: ein Element $n$ von $\mathcal{S}$ gehört dann und nur dann der Reihe $\mathcal{N}$ an, wenn $n$ Element \emph{jedes solchen} Theils $\mathcal{K}$ von $\mathcal{S}$ ist, \sina{welcher}{welches} die doppelte Eigenschaft besitzt, daß das Element $1$ in $\mathcal{K}$ enthalten ist, und daß das Bild $\phi(\mathcal{K})$ Theil von $\mathcal{K}$ ist. Zu meiner Kunstsprache: $\mathcal{N}$ ist die Gemeinheit $1_0$ oder $\phi_0(1)$ aller derjenigen Ketten $\mathcal{K}$ (in $\mathcal{S}$), \sina{in denen}{indem} das Element $1$ enthalten ist. Erst hierdurch ist der vollständige Charakter der Reihe $\mathcal{N}$ festgestellt. --- Hierzu bemerke ich beiläufig Folgendes. Die \anf{Begriffsschrift} und die \anf{Grundlagen der Arithmetik} von Frege sind zum ersten Male im letzten Sommer (1889) auf kurze Zeit in meine Hände gelangt, und ich habe mit Vergnügen gesehen, daß seine Art, das mittelbare Folgen eines Elementes auf ein anderes in einer Reihe zu erklären, im \emph{Wesentlichen} mit meinen Kettenbegriffen (37, 44) übereinstimmt; man muß sich nur durch seine etwas unbequeme Ausdrucksweise nicht zurückschrecken laßen. --- \seiteo{7}
\item Nachdem in meiner Analyse der wesentliche Charakter des einfach unendlichen Systems, deßen abstracter Typus die Zahlenreihe $\mathcal{N}$ ist, erkannt war (71, 73), fragte es sich: \emph{existirt} überhaupt ein solches System in unserer Gedankenwelt? Ohne den logischen \sina{Existenz-Beweis}{Existenzbeweis} würde es immer zweifelhaft bleiben, ob nicht der Begriff eines solchen Systems vielleicht innere Widersprüche enthält. Daher die Nothwendigkeit solcher Beweise (66, 72 meiner Schrift). 
\item Nachdem auch dies festgestellt war, fragte es sich: liegt in dem Bisherigen auch eine ausreichende \sina{\emph{Beweismethode}}{Beweismethode}, um Sätze, die für \emph{alle} Zahlen $n$ gelten sollen, allgemein zu beweisen? Ja! Der berühmte \sina{Induktions-Beweis}{Induktionsbeweis} ruht auf \sina{der}{dieser} sicheren Grundlage des \sina{Ketten-Begriffs}{Kettenbegriffes} (59, 60, 80 meiner Schrift).
\item Endlich: ist es auch möglich, die \emph{Definitionen} für \sina{Zahlen-Operationen}{Zahlen Operationen} widerspruchsfrei für \emph{alle} Zahlen $n$ aufzustellen? Ja! Dies wird durch den Satz 126 meiner Schrift in der That geleistet. ---
\end{asparaenum}\par
Damit war die Analyse beendigt, und der synthetische Aufbau konnte beginnen; \sina{er}{es} hat mir doch noch Mühe genug gemacht! Auch der \seiteo{8} Leser meiner Schrift hat es wahrlich nicht leicht; außer dem gesunden Menschenverstande gehört auch noch ein sehr starker Wille dazu, um Alles vollständig durchzuarbeiten.\par
Ich \sina{wende}{werde} mich nun noch zu einigen Stellen Ihrer Abhandlung, die ich \sina{in}{im} meiner neulichen Entgegnung nicht erwähnt habe, weil sie weniger wichtig sind; vielleicht werden aber meine darauf bezüglichen Bemerkungen noch Einiges zur Klärung der Sache beitragen.
\begin{asparaenum}[a)]
\item S.\,121, Z.\,19. Weshalb wird hier von einem \emph{Theile} gesprochen? Eine \emph{Anzahl} schreibe ich später (161 meiner Schrift) jedem \emph{endlichen} Systeme und nur einem solchen zu.
\item S.\,122, Z.\,8. Hier findet sich eine Verwechselung zwischen \emph{Abbildung} und \emph{Bild}; statt \anf{Abbildung $\overline{\phi}(\mathcal{S}')$}, müßte es heißen \anf{Abbildung $\overline{\phi}$ des Systems $\mathcal{S}'$}. Nicht $\overline{\phi}(\mathcal{S}')$, sondern $\overline{\phi}$ ist eine \emph{Abbildung} (der abbildende Maler), die aus dem \emph{System} (Original) $\mathcal{S}'$ das \sina{\emph{Bild}}{Bild} $\overline{\phi}(\mathcal{S}') = \mathcal{S}$ erzeugt. Solche Verwechselungen können aber bei unserer Untersuchung recht gefährlich werden.\seiteo{9}
\item S.\,123, Z.\,1--2. Diese Worte mögen vielleicht auf Frege paßen, auf mich gewiß nicht. Die \emph{Zahl} $1$ als Grundelement der Zahlenreihe wird von mir mit vollkommener Bestimmtheit erklärt in 71, 73, und die \emph{Anzahl} $1$ ergibt sich ebenso im Satze 164 als Folge der allgemeinen Erklärung 161. Hierzu \emph{darf} gar nichts weiter hinzugefügt werden, wenn nicht eine Trübung eintreten soll.
\item S.\,123, Z.\,29--31. Dies ist schon durch die vorhergehende Bemerkung c) erledigt. Und wie würde wohl die größere Sicherheit und die geringere Weitläufigkeit sich \emph{thatsächlich} gestalten?
\item S.\,124, Z.\,21--24. Der Sinn dieser Zeilen (sowie der vorhergehenden und folgenden) ist mir nicht ganz deutlich. Soll hier etwa der Wunsch ausgesprochen sein, meine Definition der Zahlenreihe $\mathcal{N}$ und der Aufeinanderfolge des Elementes $n'$ auf das Element $n$ womöglich anzulehnen an eine \emph{anschauliche} Reihe? Dem würde ich mich mit größter Bestimmtheit widersetzen, weil sofort die Gefahr entsteht, aus einer solchen Anschauung vielleicht unbewußt auch Sätze als selbstverständlich zu entnehmen, die vielmehr ganz abstract aus der logischen Definition von $\mathcal{N}$ abgeleitet werden müßen. Wenn ich $n'$ das auf $n$ \emph{folgende} Element \emph{nenne} (73), so soll das lediglich ein \seiteo{10} neuer \emph{Kunstausdruck} sein, durch deßen Benutzung ich nur einige Abwechselung in meine \emph{Sprache} bringe; diese Sprache würde noch einförmiger und abschreckender klingen, wenn ich auf diese Abwechselung verzichten und $n'$ immer nur das \emph{Bild} $\phi(n)$ von $n$ nennen müßte. Aber der eine Ausdruck soll genau daßelbe \emph{bedeuten} wie der andere.
\item S.\,124, Z.\,33 -- S.\,125, Z.\,7. \sina{Das}{das} in \sina{der dritten Zeile}{dem dritten Teile} meiner Erklärung 73 gewählte Wort \anf{lediglich} soll doch offenbar die einzige \emph{Einschränkung} bezeichnen, welcher das unmittelbar vorhergehende Wort \anf{gänzlich} zu unterwerfen ist; ließe man diese Einschränkung fallen, nähme also das Wort \anf{gänzlich} in seinem \emph{vollen} Sinne, so würde auch die Unterscheidbarkeit der Elemente wegfallen, welche doch für den Begriff des einfach unendlichen Systems unentbehrlich ist. Mir scheint daher dieses \anf{lediglich} durchaus nicht überflüßig, sondern nothwendig zu sein. Ich verstehe nicht, wie dies einen Anstoß erregen kann. ---
\end{asparaenum}\par
Indem ich meinen zu Anfang geäußerten Wunsch wiederhole und Sie bitte, die Ausführlichkeit meiner Erörterungen entschuldigen zu wollen, verbleibe ich mit größter Hochachtung\par\vspace{-1.5em}
\begin{center}\begin{longtable}{cp{4em}c}
Braunschweig, && Ihr \\
27.\,Februar 1890. && ergebenster \\
Petrithorpromenade 24. && R. Dedekind. \\
\end{longtable}
\end{center}\par

\def\und{\,/\,}

\end{document}